\documentclass[a4paper,12pt]{article}

\usepackage{amsthm,amsmath,stmaryrd,bbm,hyperref,geometry,color,authblk}
\usepackage[utf8]{inputenc}
\usepackage{amssymb}
\usepackage[english]{babel}
\usepackage{graphicx}
\usepackage{amsfonts,amssymb}
\usepackage{verbatim}
\usepackage{enumitem}

\newcommand{\po}{\left(}
\newcommand{\pf}{\right)}
\newcommand{\co}{\left[}
\newcommand{\cf}{\right]}

\newcommand{\R}{\mathbb R}

\newcommand{\dd}{\mathrm{d}}
\newcommand{\na}{\nabla}
\newcommand{\Id}{\mathrm{Id}}

\newtheorem{thm}{Theorem}
\newtheorem{assu}{Assumption}
\newtheorem*{assu*}{Assumption}

\newtheorem{prop}[thm]{Proposition}
\newtheorem{rem}{Remark}
\newtheorem{ex}{Example}

\title{$L^2$ geometric ergodicity for the kinetic Langevin process with non-equilibrium  steady states}
\author{Pierre Monmarché}

\begin{document}
\maketitle

\begin{abstract}
In non-equilibrium statistical physics models, the invariant measure $\mu$ of the process does not have an explicit density. In particular the adjoint $L^*$ in $L^2(\mu)$ of the generator $L$ is unknown and many classical techniques fail in this situation. An important progress has been made in \cite{HuangKopferRen}  where functional inequalities are obtained  for non-explicit steady states of kinetic equations under rather general conditions. However in \cite{HuangKopferRen}  in the kinetic case the geometric ergodicity is only deduced from the functional inequalities for the case with conservative forces, corresponding to explicit steady states. In this note we obtain $L^2$ convergence rates in the non-equilibrium case. 
\end{abstract}

\section{Motivation and result}

Consider on $\R^{d}\times \R^d$ the kinetic Langevin process $(Z_t)_{t\geqslant0}=(X_t,V_t)_{t\geqslant 0}$ solving
\begin{equation}
\label{eq:EDSLangevin}
\left\{
\begin{array}{lcl}
\dd X_t &=&  V_t \dd t \\
\dd V_t &=& b(X_t,V_t)\dd t + \sqrt{2}\sigma \dd W_t\,, 
\end{array}\right.
\end{equation}
where $b\in\mathcal C^1(\R^{2d},\R^d)$, $W$ is a $d$-dimensional Brownian motion and $\sigma \in \R_+^*$. The associated Markov generator $L$ is given by
\[L\varphi(x,v) =  v\cdot \na_x \varphi(x,v) - b(x,v)\cdot \na_v \varphi(x,v) + \sigma^2 \Delta_v \varphi(x,v)\,,\]
and the law $f_t$ of $(X_t,V_t)$  solves
\begin{equation}
\label{eq:EDP}
\partial_t f_t + v\cdot \na_x f_t = \na_v \cdot \po b f_t \pf + \sigma^2 \Delta_v f _t\,.
\end{equation}
We work under the following   condition:
\begin{assu}
\label{assu}
There  exists $\kappa,R>0$ and a positive definite symmetric matrix $A \in \R^{(2d)\times(2d)}$ such that for all $z=(x,v)$ and $z'=(x',v')$ in $\R^{2d}$ with $|z-z'|\geqslant R$,
\begin{equation}
\label{eq:condition}
\begin{pmatrix}
v-v'\\b(x,v) - b(x',v')
\end{pmatrix}\cdot A (z-z')   \leqslant - \kappa |z-z'|^2  \,.
\end{equation}
Moreover, $\|\na b\|_\infty<\infty$.
\end{assu}

\begin{ex}
For $\gamma>0$, as computed in  e.g.~\cite[Proposition 4]{M27}, there exist a positive definite $A$ and $\kappa>0$ such that for all $z=(x,v)$,
\begin{equation*}
\begin{pmatrix}
v\\ -x-\gamma v
\end{pmatrix}\cdot A z   \leqslant - \kappa |z|^2   \,.
\end{equation*}
As a consequence, Assumption~\ref{assu} holds for $b(x,v) = -x - \gamma v + F(x,v)$   as soon as $|\na F(z)|\leqslant \kappa'$ for all $z$ with $|z|\geqslant M$ for some $M>0$ and $\kappa'<\kappa/|A|$.  Indeed, in that case, we can bound, for $z,z'$ with $|z-z'|\geqslant R$,
\[|F(z)-F(z')| \leqslant |z-z'|\int_0^1 |\na F\po tz+(1-t)z'\pf| \dd t \leqslant |z-z'| \po \frac{2M}{R} \|\na F\|_\infty  + \kappa'  \pf\,.  \]
By chosing $R$ large enough so that $\varepsilon:= \kappa - |A|\po  2M \|\na F\|_\infty/R  + \kappa'\pf>0$, for $z,z'$ with $|z-z'|\geqslant R$,
\begin{equation*}
\begin{pmatrix}
v-v'\\b(x,v) - b(x',v')
\end{pmatrix}\cdot A (z-z')   \leqslant - \kappa |z-z'|^2 + |z-z'||A| |F(z)-F(z')| \leqslant  -\varepsilon|z-z'|^2   \,.
\end{equation*} 
\end{ex}

The condition~\eqref{eq:condition}  means that the norm $z\mapsto \sqrt{z\cdot A z}$ between two copies of~\eqref{eq:EDSLangevin} driven by the same Brownian motion deterministically decays at constant rate as long as the processes are at a distance larger than $R$. When $R=0$ the system is globally dissipative (or contracting) and the question of the long-time convergence of the process is well-understood \cite{M27}.  

 More generally, under Assumption~\ref{assu}, the process~\eqref{eq:EDSLangevin}  is known to be  ergodic and thus admits a unique invariant measure $\mu$, which is thus the unique steady state of~\eqref{eq:EDP}. The operator $L$ being hypoelliptic and the process being controllable, $\mu$ admists a positive density (still denoted by $\mu$).   In the so-called conservative (or equilibrium) case where $b(x,v) = -\na U(x) -\gamma v$ for some $U\in\mathcal C^2(\R^d,\R)$ and $\gamma>0$, then $\mu$ is the Gibbs measure with density proportional to $\exp(-\frac{\gamma}{\sigma^2} H)$, with $H(x,v) = U(x) + |v|^2/2$. In particular, under $\mu$, the position $X$ and the velocity are independent. This is not true in general in the non-equilibrium case where $b$ is not of the previous form \cite{Mouad}. In this situation in general the measure $\mu$ has no explicit density. For motivations on this non-equilibrium situation, we refer to \cite{Iacobucci,iacobucci2021thermo,CuneoEckmann} and references within.

Under Assumption~\ref{assu}, with a probabilistic approach, the long-time convergence of $\nu_t$ to $\mu$ is known to be exponentially fast in $V$-norms using Harris theorem \cite{talay2002stochastic,MATTINGLY2002185} or in Wasserstein 1 distance using reflection couplings \cite{Eberle}.  Alternatively, with a PDE point of view, a natural way to quantify the convergence of $\nu_t$ to $\mu$ is the $L^2(\mu)$ norm of the relative density $h_t=f_t/\mu$, 
\[\|h_t - 1\|_{2}^2 = \int_{\R^{2d}} \po h_t - 1 \pf^2 \mu = \int_{\R^{2d}} \frac{(f_t - \mu)^2}{\mu}\,.\]
The evolution of $h_t$ is governed by the  dual of $L$ in $L^2(\mu)$, 
\[L^* \varphi =  - v\cdot \na_x \varphi + \po -  b + 2 \sigma^2 \na_v \ln \mu\pf \cdot \na_v \varphi + \sigma^2 \Delta_v \varphi\,,\]
since
\[\partial_t h_t = L^* h_t\,.\]
However, when $\mu$ is not explicit, neither is $L^*$. 

Apart the globally dissipative case, or for small perturbations of an equilibrium situation as in \cite{Iacobucci}, the recent  works \cite{dietert20232,MonmarcheWang} establish an exponential convergence of $\|h_t - 1\|_{2}^2$ in non-equilibrium situations under some  restrictive conditions. In the present note we prove it under the natural Assumption~\ref{assu} (which is not covered by \cite{dietert20232,MonmarcheWang}). The result is obtained by combining the functionnal inequalities established in the very recent~\cite{HuangKopferRen} with the modified norm result of \cite{Gamma}. The work~\cite{HuangKopferRen} is concerned with more general settings with possibly non-constant diffusion matrix and a small non-linear perturbation where the drift and diffusion matrix can slightly depend on the law of the process itself. However, in the kinetic case, in terms of long-time convergence, only the globally dissipative and equilibrium cases are covered, in \cite[Theorems 3.7 and 3.9]{HuangKopferRen}. This is because, in the non-globally dissipative case, the proof is based on the modified norm result of \cite{villani2009hypocoercivity}, which requires $\mu$ and $L^*$ to be explicit.

An immediate corollary of  \cite[Theorem 3.5]{HuangKopferRen} is the following:

\begin{prop}\label{prop}
Under Assumption~\ref{assu}, $\mu$ satisfies a Poincaré inequality, namely there exists $C>0$ such that for all $g\in \mathcal H^1(\mu)$ with $\int_{\R^{2d}} g  \mu=0$,
\begin{equation}
\label{eq:Poincare}
\|g\|_2^2 \leqslant C \|\na g\|_2^2\,.
\end{equation}
\end{prop}

\begin{rem}\label{rem}
Notice that~\eqref{eq:Poincare} is a classical Poincaré inequality with respect to the full gradient $|\na g|^2$, and not a Poincaré inequality associated to the carré du champ of the process (as in \cite{BakryGentilLedoux}), which for \eqref{eq:EDSLangevin} is $\Gamma(g) = |\na_v g|^2$ (for which no such inequality holds, as can be seen by considering a non-constant $g$ depending only on $x$).
\end{rem}

In fact, under more general settings, \cite[Theorem 3.5]{HuangKopferRen} shows that $\mu$ satisfies a log-Sobolev inequality, which is stronger than the Poincaré inequality (see e.g. \cite[Proposition 5.1.3]{BakryGentilLedoux}). Using this log-Sobolev inequality with \cite[Theorem 10]{Gamma} gives an exponential convergence in relative entropy for the semi-group $P_t = e^{tL}$, i.e. $P_t \varphi(x) = \mathbb E_x (\varphi(X_t))$.  However it is unclear how to obtain a convergence in relative entropy for the law $\nu_t$ itself, which would be the question of interest. This is why we stick to the $L^2$ framework where the operator norms $\|e^{tL}-\mu\|_2$ and $\|e^{tL^*}-\mu\|_2$ are the same by duality. This yields the following:

\begin{thm}\label{thm}
Under Assumption~\ref{assu}, there exists $c>0$ such that 
\[\|h_t - 1 \|_{2} \leqslant e^{-c \min(t,t^3) } \|h_0 - 1 \|_{2} \,.\]
\end{thm}

In relation with Remark~\ref{rem}, let us discuss now that for any diffusion process with Lipschitz coefficient and bounded diffusion matrix, a decay of $\|P_t - \mu \|_2$ implies a Poincaré inequality with respect to the full gradient (independently from the carré du champ of the process). Consider on $\R^d$ a non-explosive diffusion process solving
\begin{equation}
\label{eq:Z}
\dd Z_t = b(Z_t) \dd t + \sigma(Z_t) \dd B_t\,,
\end{equation}
where $b\in\mathcal C^1(\R^d,\R^d)$ and $\sigma\in\mathcal C^1(\R^d,\R^{d\times d})$. Denote by $(P_t)_{t\geqslant 0}$ the associated Markov semi-group. 

\begin{prop}\label{prop2}
Assume that the process \eqref{eq:Z} admits an invariant measure $\mu$,  $\sigma$ is bounded, there exists $K>0$ such that
\begin{equation}
\label{eq:onesided}
(x-y)\cdot \po b(x)-b(y)\pf  + \|\sigma(x)-\sigma(y)\|_{HS}^2 \leqslant K |x-y|^2 \qquad \forall x,y\in\R^d
\end{equation}
and  there exists $t_0>0$ such that $\|P_{t_0} - \mu\|_2 < 1$. Then $\mu$ satisfies a Poincaré inequality~\eqref{eq:Poincare} with constant
\[C = \|\sigma\|_\infty^2 \frac{e^{2K t_0}-1}{2K(1-\|P_{t_0} - \mu\|_2^2)}  \,.  \] 
\end{prop}

\section{Proofs}\label{sec:proof}

\begin{proof}[Proof of Theorem~\ref{thm}]
Denoting by $P_t^*$ the dual of $P_t$ in $L^2(\mu)$,
\[\|h_t - 1\|_{2} \leqslant \|P_t^*- \mu\|_{2} \|h_0-1\|_2 = \|P_t- \mu\|_{2}^2 \|h_0-1\|_2\,. \]
The result is thus a direct consequence of \cite[Theorem 10]{Gamma}, which can be applied thanks to Proposition~\ref{prop} and the condition that $\|\na b\|_\infty<\infty$. Since \cite[Theorem 10]{Gamma} is a general statement, for the reader's convenience we provide here the main lines  of proof in the specific case~\eqref{eq:EDSLangevin}. Fix $g \in \mathcal H^1(\mu)$ smooth and set $g_t = P_t g - \int g\mu$ for $t\geqslant 0$. Using that $\int_{\R^{2d}} L \varphi \mu=0$ for all smooth $\varphi$,
\[\partial_t \|g_t\|_2^2 = 2 \int_{\R^{2d}} g_t L g_t \mu =  \int_{\R^{2d}} \po  -L(g_t^2) + 2 g_t L g_t \pf  \mu = - 2\sigma^2 \int_{\R^{2d}} |\na_v g_t|^2 \mu\,.\]
Similarly, writing
\[J(x,v) = \begin{pmatrix}
0 & \na_x b(x,v) \\
\Id & \na_v b(x,v) 
\end{pmatrix}\]
the Jacobian matrix of the drift of~\eqref{eq:EDSLangevin}, for any constant matrix $D\in\R^{(2d)\times(2d)}$,
\begin{align*}
\partial_t \| D\na g_t\|_2^2 &= 2 \int_{\R^{2d}} D\na g_t \cdot D  \na L g_t \\
&=   \int_{\R^{2d}} \po - L(|D\na g_t|^2) + 2  D\na g_t  \cdot  L  D   \na g_t + 2D\na g_t \cdot  D J\na g_t   \pf \mu  
\end{align*}
where, for a vector-valued function $m=(m_1,\dots,m_n)$, $Lm$ means $(Lm_1,\dots,L m_n)$, and relying on the fact that $D$ and $\sigma$ are constant, we used that
\[D \na L g_t = LD\na g_t + D J \na g_t\,.\]
As a consequence,
\[ \partial_t \| D\na g_t\|_2^2  = - 2 \sigma^2 \sum_{i=1}^d \|\partial_{v_i} D\na g_t\|_2^2 +  2\int_{\R^{2d}}   D\na g_t \cdot  D J\na g_t \mu \leqslant 2\int_{\R^{2d}}   D\na g_t \cdot D J\na g_t \mu\,. \]  
Set $\alpha(t) = 1-e^{-t/3}$  and take $D_t$ as the symmetric square-root of
\[D_t^2 = \varepsilon \begin{pmatrix}
\alpha^3(t) \Id & -\alpha^2(t) \Id \\
-\alpha^2(t) \Id & \alpha(t) \Id   
\end{pmatrix}\,,\]
for some $\varepsilon>0$ to be fixed later on. Write
\[\mathcal N_t = \|g_t\|_2^2 + \| D_t \na g_t\|_2^2 = \|g_t\|_2^2 + \alpha(t) \| (\na_v - \alpha (t)\na_x)  g_t\|_2^2\,.\]
Gathering the previous computations,
\[\partial_t\mathcal N_t \leqslant \int_{\R^d} \na g_t R_t \na g_t \,, \]
with, writing $\alpha$ instead of $\alpha(t)$ to alleviate notations,
\begin{align*}
R_t & = -2\sigma^2\begin{pmatrix}
0 & 0 \\ 0 & \Id 
\end{pmatrix} + \partial_t(D_t^2) + 2 D_t^2 J \\
&= \varepsilon \begin{pmatrix}
(3\alpha'-2)\alpha^2 \Id  & 2 \alpha^3 \na_x b - 2\alpha^2 \na_v b -2\alpha'\alpha \Id  \\ 2(1-\alpha') \alpha \Id  & -\alpha^2 \na_x b + \alpha \na_v b + \po \alpha' -2\sigma^2\varepsilon^{-1}\pf  \Id 
\end{pmatrix}
\end{align*}
Using that  $\alpha' \in(0, 1/3]$ and $\alpha\leqslant 1$, for any $z=(x,v)\in\R^{2d}$,
\begin{align*}
z\cdot R_t z &\leqslant - \alpha^2 \varepsilon |x|^2 + 2\po  \|\na_x b\|_\infty +\|\na_v b\|_\infty +  1\pf \alpha \varepsilon |x||v| + \co \varepsilon\po \|\na_x b\|_\infty +\|\na_v b\|_\infty + 1  \pf - 2\sigma^2 \cf |v|^2 \\
& \leqslant - \frac12 \alpha^2 \varepsilon |x|^2 + (M\varepsilon -2\sigma^2) |v|^2 
\end{align*}
for some $M>0$ independent from $t$ and $\varepsilon$. We take $\varepsilon$ small enough to get
\[\partial_t \mathcal N_t \leqslant -\frac{\varepsilon\alpha^2}{2} \|\na g_t\|_2^2 \leqslant - \frac{\varepsilon\alpha^2}{2C+4\varepsilon} \mathcal N_t, \] 
where we used in the second inequality the Poincaré inequality from Proposition~\ref{prop} and that $|D_t|^2 \leqslant 2 \varepsilon$ for all $t\geqslant 0$.
As a conclusion,
\[\|g_t\|_2^2 \leqslant \mathcal N_t \leqslant \exp\po - \frac{\varepsilon}{2C+4\varepsilon} \int_0^t \alpha^2(s)\dd s\pf \mathcal N_0 \leqslant e^{-c\min(t,t^3)} \|g_0\|_2^2\,, \]
for some $c>0$, where we used that $D_0=0$. 
 \end{proof}

  \begin{proof}[Proof of Proposition~\ref{prop2}]
  Thanks to~\eqref{eq:onesided}, along a synchronous coupling of two processes, 
  \[\mathbb E\po |Z_t - Z_t'|^2\pf \leqslant e^{2 K t} \mathbb E\po |Z_0 - Z_0'|^2\pf\]
 see e.g. \cite[Theorem 2.5]{FYWang}. By \cite{Kuwada1} and the invariance of $\mu$ by $P_t$, this implies that
  \[\|\na P_t \varphi\|^2_2 \leqslant e^{2K t} \|\na \varphi\|_2\]
  for all $\varphi\in H^1(\mu)$ and $t\geqslant 0$.  Then, writing $g_t = P_t g - \mu(g)$ for some $g\in H^1(\mu)$,
  \[ \|g_0\|_2^2 - \|g_t\|_2^2 = 2 \int_0^t \|\sigma\na g_s\|_2^2 \dd s  \leqslant 2\|\sigma\|_\infty \frac{e^{2Kt}-1}{2K} \|\na g\|_2^2\,.\]
Applying this at time $t_0$,  conclusion then follows from
  \[\|g_{t_0}\|_2^2 \leqslant \|P_{t_0}-\mu\|_2^2 \|g_0\|_2^2\,.\] 
  \end{proof}

\subsection*{Acknowledgements}

P.M. would like to thank Feng-Yu Wang for pointing him out the reference~\cite{HuangKopferRen}. His research is supported by the projects SWIDIMS (ANR-20-CE40-0022) and CONVIVIALITY (ANR-23-CE40-0003) of the French National Research Agency and by the project EMC2 from the European Union’s Horizon 2020 research and innovation program (grant agreement No 810367).
  
\bibliographystyle{plain} 
\bibliography{biblio}

\end{document}